\documentclass[12pt]{article}
\usepackage{a4}
\usepackage{amsfonts}
\usepackage{amssymb}

\begin{document}

\title{\textbf{On the smallest poles of topological zeta functions}}
\author{Dirk Segers \and Willem Veys}

\maketitle

\begin{abstract}
We study the local topological zeta function associated to a
complex function that is holomorphic at the origin of
$\mathbb{C}^2$ (respectively $\mathbb{C}^3$). We determine all
possible poles less than $-1/2$ (respectively $-1$). On
$\mathbb{C}^2$ our result is a generalization of the fact that the
log canonical threshold is never in $]5/6,1[$. Similar statements
are true for the motivic zeta function.
\end{abstract}

\section{Introduction}
\noindent \textbf{(1.1)} Let $f$ be the germ of a holomorphic
function on a neighbourhood of the origin $0$ in $\mathbb{C}^n$
which satisfies $f(0)=0$ and which is not identically zero. Let
$g:V \rightarrow U \subset \mathbb{C}^n$ be an embedded resolution
of a representative of $f^{-1}\{0\}$. We denote by $E_i$, $i \in
T$, the irreducible components of $g^{-1}(f^{-1}\{0\})$, and by
$N_i$ and $\nu_i-1$ the multiplicities of $f \circ g$ and
$g^*(dx_1 \wedge \cdots \wedge dx_n)$ along $E_i$. The
$(N_i,\nu_i)$, $i \in T$, are called the numerical data of the
resolution $(V,g)$. For $I \subset T$ denote also $E_I:=\cap_{i
\in I} E_i$ and $\stackrel{\circ}{E_I} := E_I \setminus (\cup_{j
\notin I} E_j)$.

The set of germs of holomorphic functions on a neighbourhood of $0
\in \mathbb{C}^n$ will be denoted by $\mathcal{O}_n$.

\vspace{0,5cm}

\noindent \textbf{(1.2)} To $f$ one associates the local
topological zeta function
\[ Z_f(s) = Z_{\mathrm{top},0,f}(s) := \sum_{I \subset T}
\chi(\stackrel{\circ}{E_I} \cap g^{-1}\{0\}) \prod_{i \in I}
\frac{1}{\nu_i+sN_i}.\] Here $s$ is a complex variable and
$\chi(\cdot)$ denotes the topological Euler-Poincar\'{e}
characteristic. The remarkable fact that $Z_f(s)$ does not depend
on the chosen resolution was first proved in \cite{DenefLoeser1}
by expressing it as a limit of Igusa's $p$-adic zeta functions.

\vspace{0,5cm} \footnoterule{\noindent \footnotesize{2000
\emph{Mathematics Subject Classification.} Primary 14B05 14J17
32S05; Secondary 14E15 14H20} \\ \emph{Key words.} Topological
zeta function, resolution of singularities, log canonical
threshold.}
\newpage

\noindent \textbf{(1.3)} The log canonical threshold $c_0(f)$ of
$f$ at $0 \in \mathbb{C}^n$ is by definition
\[ \sup \{ c \in \mathbb{Q} \mid \mbox{ the pair } (\mathbb{C}^n , c \; \mbox{div} \; f)
\mbox{ is log canonical in a neighbourhood of } 0 \}. \] We can
describe it (see \cite[Prop 8.5]{Kollar2}) in terms of the
embedded resolution $(V,g)$ as $c_0(f)=\min\{ \nu_i/N_i \mid i \in
T \}$. In particular, this minimum is independent of the chosen
resolution. Consequently, $-c_0(f)$ is the largest candidate pole
of $Z_f(s)$. The log canonical threshold has already been studied
in various papers of Alexeev, Ein, Koll\'ar, Kuwata, Musta\c t\u
a, Prokhorov, Reid, Shokurov and others; especially the sets \[
\mathcal{T}_n := \{ c_0(f) \mid f \in \mathcal{O}_n \}, \] with $n
\in \mathbb{Z}_{>0}$, are the subject of interesting conjectures.

It is natural to investigate whether more quotients $-\nu_i/N_i$,
$i \in T$, yield invariants of the germ of $f$ at $0$. Of course,
the whole set $\{ -\nu_i/N_i \mid i \in T \}$ depends on the
chosen resolution (for n=2 however one could consider such a set
associated to the minimal resolution); but its subset consisting
of the poles of $Z_f(s)$ is an invariant of $f$. Philosophically,
these poles are induced by `important' components $E_i$, which
occur in every resolution. For $n \in \mathbb{Z}_{>0}$, we define
the set $\mathcal{P}_n$ by
\[ \mathcal{P}_n := \{ s_0 \mid \exists f \in \mathcal{O}_n \,
: \, Z_{f}(s) \textsl{\mbox{ has a pole in }} s_0 \}. \] The case
$n=1$ is trivial: $\mathcal{T}_1=\{1/i \mid i \in
\mathbb{Z}_{>0}\}$ and $\mathcal{P}_1= \{-1/i \mid i \in
\mathbb{Z}_{>0}\}$.

\vspace{0,5cm}

\noindent \textbf{(1.4)} When $n=2$, it is known that
$\mathcal{T}_2 \cap ]5/6,1[ = \emptyset$ (see \cite{Reid}). Because it follows from \cite{Veysdetermination}
that $-c_0(f)$ is a pole (and thus the largest pole) of $Z_f(s)$,
the statement $\mathcal{P}_2 \cap ]-1,-5/6[ = \emptyset$ would be
a remarkable generalization. It is in fact not hard to prove this
generalization. In this article, we will prove more:
\begin{eqnarray} \label{result1}
\mathcal{P}_2 \cap ]-\infty,-1/2[ & = & \{-1/2-1/i \mid i \in
\mathbb{Z}_{>1} \} \\  & = & \{-1,-5/6,-3/4,-7/10,\ldots \}.
\nonumber
\end{eqnarray}

\vspace{0,5cm}

\noindent \textbf{(1.5)} Koll\'ar proved in \cite{Kollar1} that
$\mathcal{T}_3 \cap ]41/42,1[ = \emptyset$. It turns out that
there is no analogous result for $\mathcal{P}_3$. Actually, we
will give examples of zeta functions with  poles in $]-1,-41/42[$
which are moreover arbitrarily near to $-1$. On the other hand, we
prove the analogue of (\ref{result1}), which appears to be
\begin{eqnarray}
\mathcal{P}_3 \cap ]-\infty,-1[ & = & \{-1-1/i \mid i \in
\mathbb{Z}_{>1} \}.
\end{eqnarray}
In general, we expect that $\mathcal{P}_n \cap ]-\infty,-(n-1)/2[
= \{-(n-1)/2-1/i \mid i \in \mathbb{Z}_{>1} \}$.

\vspace{0,2cm} \noindent \textsl{Remark.} One can easily show that
$\mathcal{P}_n \cap ]-\infty,-n+1[ = \emptyset$ if $n \geq 2$.

\section{Curves}
\noindent \textbf{(2.1)} We will determine $\mathcal{P}_2 \cap
]-\infty,-1/2[$.  Let $f$ be the germ of a holomorphic function on
a neighbourhood of the origin $0$ in $\mathbb{C}^2$ which
satisfies $f(0)=0$ and which is not identically zero. Let $(V,g)$
be the \textsl{minimal} embedded resolution of $f^{-1}\{0\}$.
Write $g=g_1 \circ \cdots \circ g_t$ as a composition of
blowing-ups $g_i$, $i \in T_e:=\{1,\ldots,t\}$. The exceptional
curve of $g_i$ and also the strict transforms of this curve are
denoted by $E_i$. The irreducible components of $f^{-1}\{0\}$ and
their strict transforms are denoted by $E_j$, $j \in T_s$.

\vspace{0,5cm}

\noindent \textbf{(2.2)} The dual (minimal) embedded resolution
graph of $f^{-1}\{0\}$ is obtained as follows. One associates a
vertex to each exceptional curve in the minimal embedded
resolution (represented by a dot), and to each branch of the
strict transform of $f^{-1}\{0\}$ (represented by a circle). One
also associates to each intersection an edge, connecting the
corresponding vertices. The fact that $E_i$ has numerical data
$(N_i,\nu_i)$ is denoted by $E_i(N_i,\nu_i)$.

\vspace{0,5cm}

\noindent \textbf{(2.3)} Let $E_i$ be an exceptional curve and let
$E_j$, $j \in J$, be the components that intersect $E_i$ in $V$.
Set $\alpha_j=\nu_j-(\nu_i/N_i)N_j$ for $j \in J$. Then we have
the relation
\begin{equation}\label{relatiealpha}
\sum_{j \in J} (\alpha_j-1)+2=0,
\end{equation}
which was first proved by Loeser in \cite{Loeser}, and later more
conceptually by the second author in \cite{Veysrelations}.
\\ Suppose that $\alpha_j \not= 0$, which is equivalent to
$-\nu_i/N_i \not= -\nu_j/N_j$, for all $j \in J$. Then one
computes easily that the contribution of $E_i$ to the residue
$\mathcal{R}$ of $Z_f(s)$ at the candidate pole $-\nu_i/N_i$ is
\begin{equation}\label{contributie}
\frac{1}{N_i} \left( \chi(\stackrel{\circ}{E_{\{i\}}}) + \sum_{j
\in J} \alpha_j^{-1} \right)
\end{equation}
(see \cite[section 2.3]{Veysdetermination}). From
(\ref{relatiealpha}) and (\ref{contributie}) it follows that
$\mathcal{R}=0$ if $J$ contains one or two elements. This is the
easy part of the following theorem. The other part is more
difficult and is proved in \cite{Veysdetermination}.

\vspace{0,5cm}

\noindent \textbf{(2.4) Theorem.} \textsl{We have that $s_0$ is a
pole of $Z_f(s)$ if and only if $s_0=-\nu_i/N_i$ for some
exceptional curve $E_i$ intersecting at least three times other
components, or $s_0=-1/N_j$ for some irreducible component $E_j$
of the strict transform of $f^{-1}\{0\}$.}

\vspace{0,5cm}

\noindent The following lemma is obtained by elementary
calculations.

\vspace{0,5cm}

\noindent \textbf{(2.5) Lemma.} \textsl{Suppose that we have blown
up $k$ times but we have not yet an embedded resolution. Let $P$
be a point of the strict transform of $f^{-1}\{0\}$ with
multiplicity $\mu$ in which we do not have normal crossings yet.
Let $g_{k+1}$ be the blowing-up at $P$.}

\noindent \hspace{0,5cm} \textsl{(a) Suppose that two exceptional
curves $E_i$ and $E_j$ contain $P$. Then the new candidate pole
$-\nu_{k+1}/N_{k+1}=-(\nu_i+\nu_j)/(N_i+N_j+\mu)$ is larger than
$\min\{-\nu_i/N_i,-\nu_j/N_j\}$.}

\noindent \hspace{0,5cm} \textsl{(b) Suppose that exactly one
exceptional curve $E_i$ contains $P$ and that $\mu \geq 2$. Then
$E_{k+1}$ has numerical data $(N_i+\mu,\nu_i+1)$ and
$-(\nu_i+1)/(N_i+\mu)$ is in between $-1/\mu$ and $-\nu_i/N_i$.}

\noindent \hspace{0,5cm} \textsl{(c) Suppose that exactly one
exceptional curve $E_i$ contains $P$ and that $\mu=1$. Remark that
the two curves are tangent at $P$ because we do not have normal
crossings at $P$. Let $g_{k+2}$ be the blowing-up at $E_i \cap
E_{k+1}$. Because the strict transform of $f^{-1}\{0\}$ does not
intersect $E_{k+1}$ after this blowing-up, we do not have to blow
up at a point of $E_{k+1}$ anymore. Because $E_{k+1}$ is
intersected once, it follows from (2.3) that the contribution of
$E_{k+1}$ to the residue at the candidate pole
$-\nu_{k+1}/N_{k+1}$ is zero. The numerical data of $E_{k+2}$ are
$(2N_i+2,2\nu_i+1)$, and $-(2\nu_i+1)/(2N_i+2)$ is in between
$-1/2$ and $-\nu_i/N_i$.}

\vspace{0,5cm}

\noindent \textbf{(2.6)} Suppose that after some blowing-ups, we
do not have normal crossings at a point $P$. Suppose also that the
candidate poles associated to the exceptional curves through $P$
are all larger than or equal to $-1/2$. Then it follows from the
above lemma that the components above $P$ in the final resolution
do not give a contribution to a pole less than $-1/2$.

\vspace{0,2cm} \noindent \textbf{Corollary.} \textsl{Zeta
functions of singularities of multiplicity at least $4$ do not
have a pole in} $]-\infty,-1/2[ \setminus \{-1\}$.

\vspace{0,2cm} \noindent Indeed, every exceptional curve in the
minimal embedded resolution of $f^{-1}\{0\}$ lies above a point of
$E_1$ (considered in the stage when it is created), which has a
candidate pole larger than or equal to $-1/2$.

\vspace{0,5cm}

\noindent \textbf{(2.7)} If $f \in \mathcal{O}_2$ has multiplicity
$2$ or $3$, we will use the Weierstrass Preparation Theorem and
coordinate transformations to obtain an `easier' element of
$\mathcal{O}_2$ with the same zeta function.

We illustrate this in the case that $f \in \mathcal{O}_2$ has
multiplicity $3$ and the homogeneous part of degree $3$ of $f$ is
$f_3=y^3+xy^2=y^2(y+x)$. According to the Weierstrass Preparation
Theorem, we have that $f=(y^3+a_1(x)y^2+a_2(x)y+a_3(x))h(x,y)$,
with $\mbox{mult}(a_1(x))=1$, $\mbox{mult}(a_2(x)) \geq 3$,
$\mbox{mult}(a_3(x)) \geq 4$ and $h(0,0) \not= 0$. Because $h(0,0)
\not= 0$, the resolutions and the local topological zeta functions
of $f$ and $y^3+a_1(x)y^2+a_2(x)y+a_3(x)$ are the same. One can
check that there exists a coordinate transformation $(x,y) \mapsto
(x,y-k(x))$ such that the last function becomes of the form
$y^3+b_1(x)y^2+b_3(x)$, with $\mbox{mult}(b_1(x))=1$ and
$\mbox{mult}(b_3(x)) \geq 4$. After another coordinate
transformation, we get the form $y^3+xy^2+g(x)$, with
$\mbox{mult}(g(x)) \geq 4$.

\vspace{0,5cm}

\noindent \textbf{(2.8) Theorem.} \textsl{We have}
\[ \mathcal{P}_2 \cap \left]-\infty,-\frac{1}{2}\right[ = \left\{ \left.
-\frac{1}{2}-\frac{1}{i} \right| i \in \mathbb{Z}_{>1} \right\} \]
\textsl{and every local topological zeta function has at most one
pole in $]-1,-1/2]$.}

\vspace{0,2cm} \noindent \emph{Proof.} (a) Suppose that
$\mbox{mult}(f)$, the multiplicity of $f$ at the origin of
$\mathbb{C}^2$, is equal to $2$. Then $f$ is holomorphically
equivalent to $y^2$ or $y^2+x^k$ for some $k \in \mathbb{Z}_{>1}$.
If it is $y^2$, the only pole of $Z_f(s)$ is $-1/2$. If $k=2$, the
only pole of $Z_f(s)$ is $-1$. If $k$ is odd, write $k=2r+1$.
After $r$ blowing-ups, the strict transform of $f^{-1}\{0\}$ is
nonsingular and tangent to $E_r$. The numerical data of $E_i$,
$i=1,\ldots,r$, are $(2i,i+1)$. To get the minimal embedded
resolution, we now blow up twice. The dual resolution graph and
the numerical data are given below. \\
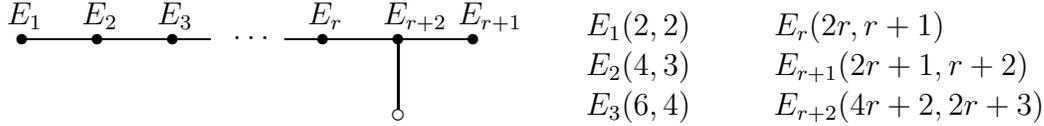
\begin{picture}(150,20)(-5,2)
\put(0,15){\line(1,0){25}} \put(28,15){$\ldots$}
\put(35,15){\line(1,0){25}} \put(0,15){\circle*{1.5}}
\put(10,15){\circle*{1.5}} \put(20,15){\circle*{1.5}}
\put(40,15){\circle*{1.5}} \put(50,15){\circle*{1.5}}
\put(60,15){\circle*{1.5}} \put(50,14.5){\line(0,-1){8.75}}
\put(50,5){\circle{1.5}} \put(-2,17){$E_1$} \put(8,17){$E_2$}
\put(18,17){$E_3$} \put(38,17){$E_r$} \put(48,17){$E_{r+2}$}
\put(58,17){$E_{r+1}$} \put(75,5){\shortstack[l]{$E_1(2,2)$ \\
$E_2(4,3)$ \\ $E_3(6,4)$}}
\put(100,5){\shortstack[l]{$E_r(2r,r+1)$
\\ $E_{r+1}(2r+1,r+2)$ \\ $E_{r+2}(4r+2,2r+3)$}}
\end{picture}
If $k$ is even and larger than $2$, write $k=2r$. Easy
calculations give the following dual resolution graph. \\
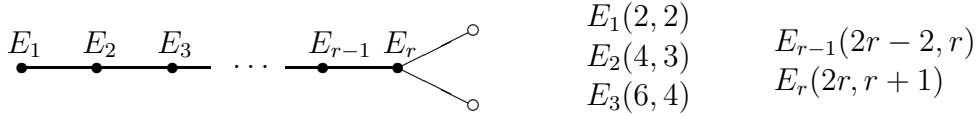
\begin{picture}(150,20)(-5,2)
\put(0,10){\line(1,0){25}} \put(28,10){$\ldots$}
\put(35,10){\line(1,0){15}} \put(0,10){\circle*{1.5}}
\put(10,10){\circle*{1.5}} \put(20,10){\circle*{1.5}}
\put(40,10){\circle*{1.5}} \put(50,10){\circle*{1.5}}
\put(50,10){\line(2,1){9.3}} \put(50,10){\line(2,-1){9.3}}
\put(60,15){\circle{1.5}} \put(60,5){\circle{1.5}}
\put(-2,12){$E_1$} \put(8,12){$E_2$} \put(18,12){$E_3$}
\put(38,12){$E_{r-1}$} \put(48,12){$E_r$}
\put(75,5){\shortstack[l]{$E_1(2,2)$ \\ $E_2(4,3)$ \\ $E_3(6,4)$}}
\put(100,7){\shortstack[l]{$E_{r-1}(2r-2,r)$
\\ $E_r(2r,r+1)$}}
\end{picture}
Because $-(2r+3)/(4r+2)=-1/2-1/(2r+1)$ and
$-(r+1)/(2r)=-1/2-1/(2r)$, it follows from (2.4) that
\begin{eqnarray*} \{ s_0 \mid \exists f \in \mathcal{O}_2 & : &
\mbox{mult}(f)=2 \mbox{ and } Z_f(s) \mbox{ has a pole in } s_0 \}
\\ & = & \left\{ \left. -\frac{1}{2}-\frac{1}{i} \right| i \in
\mathbb{Z}_{>1} \right\} \cup \left\{ -\frac{1}{2} \right\}.
\end{eqnarray*}
Remark that Newton polyhedra could also be used to deal with (a),
see \cite{DenefLoeser1}.

\vspace{0,4cm} \noindent \hspace{1cm} (b) Suppose that
$\mbox{mult}(f)=3$. Up to an affine coordinate transformation,
there are three cases for $f_3$.

\vspace{0,3cm} \noindent \hspace{0,5cm} (b.1) Case $f_3=xy(x+y)$.
After one blowing-up we get an embedded resolution. The poles of
$Z_f(s)$ are $-1$ and $-2/3=-1/2-1/6$.

\vspace{0,3cm} \noindent \hspace{0,5cm} (b.2) Case $f_3=y^2(y+x)$.
According to (2.7), we may suppose that $f=y^3+xy^2+g(x)$, where
$g(x)$ is a holomorphic function in the variable $x$ of
multiplicity $k \geq 4$. If $g(x)=0$, the poles of $Z_f(s)$ are
$-1$ and $-1/2$. Consider now the case that $k$ is odd. Write
$k=2r+1$. After $r$ blowing-ups we get an embedded resolution with
the following dual resolution graph and numerical data. \\
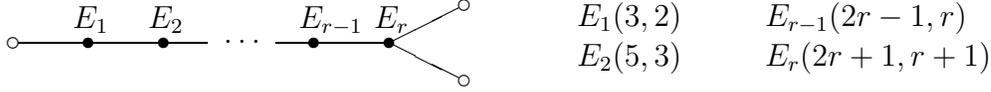
\begin{picture}(150,20)(-5,2)
\put(0.75,10){\line(1,0){25}} \put(28,10){$\ldots$}
\put(35,10){\line(1,0){15}} \put(0,10){\circle{1.5}}
\put(10,10){\circle*{1.5}} \put(20,10){\circle*{1.5}}
\put(40,10){\circle*{1.5}} \put(50,10){\circle*{1.5}}
\put(50,10){\line(2,1){9.3}} \put(50,10){\line(2,-1){9.3}}
\put(60,15){\circle{1.5}} \put(60,5){\circle{1.5}}
\put(8,12){$E_1$} \put(18,12){$E_2$} \put(38,12){$E_{r-1}$}
\put(48,12){$E_r$} \put(75,7){\shortstack[l]{$E_1(3,2)$ \\
$E_2(5,3)$}} \put(100,7){\shortstack[l]{$E_{r-1}(2r-1,r)$
\\ $E_r(2r+1,r+1)$}}
\end{picture}
If $k$ is even, write $k=2r$. After $r+1$ blowing-ups we get the
following picture.
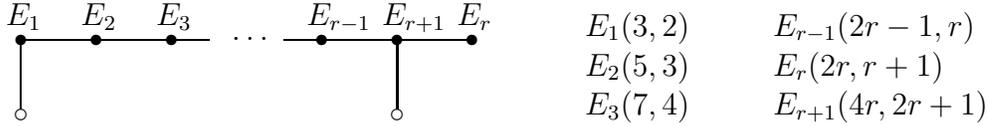
\begin{picture}(150,20)(-5,2)
\put(0,15){\line(1,0){25}} \put(28,15){$\ldots$}
\put(35,15){\line(1,0){25}} \put(0,15){\circle*{1.5}}
\put(10,15){\circle*{1.5}} \put(20,15){\circle*{1.5}}
\put(40,15){\circle*{1.5}} \put(50,15){\circle*{1.5}}
\put(60,15){\circle*{1.5}} \put(0,14.5){\line(0,-1){8.75}}
\put(50,14.5){\line(0,-1){8.75}} \put(0,5){\circle{1.5}}
\put(50,5){\circle{1.5}} \put(-2,17){$E_1$} \put(8,17){$E_2$}
\put(18,17){$E_3$} \put(38,17){$E_{r-1}$} \put(48,17){$E_{r+1}$}
\put(58,17){$E_{r}$} \put(75,5){\shortstack[l]{$E_1(3,2)$ \\
$E_2(5,3)$ \\ $E_3(7,4)$}}
\put(100,5){\shortstack[l]{$E_{r-1}(2r-1,r)$
\\ $E_r(2r,r+1)$ \\ $E_{r+1}(4r,2r+1)$}}
\end{picture}
The poles appearing in (b.2) are in the desired set because
$-(r+1)/(2r+1)=-1/2-1/(4r+2)$ and $-(2r+1)/(4r)=-1/2-1/(4r)$.

\vspace{0,3cm} \noindent \hspace{0,5cm} (b.3) Case $f_3=y^3$. We
may suppose that $f$ is of the form
\[y^3+a_4x^4+b_3yx^3+a_5x^5+b_4yx^4 +a_6x^6+b_5yx^5+\cdots,\]
where $a_i,b_i \in \mathbb{C}$. If $f=f_3=y^3$ then the only pole
of $Z_f(s)$ is $-1/3$. Otherwise there is an integer $r \geq 1$
such that after blowing up $r$ times and always taking the charts
determined by $g_i(x,y)=(x,xy)$, we get $(g_1 \circ \cdots \circ
g_r)^* dx \wedge dy = x^r dx \wedge dy$ and $f \circ g_1 \circ
\cdots \circ g_r = x^{3r}(y^3+a_{3r+1}x
+b_{2r+1}yx+a_{3r+2}x^2+b_{2r+2}yx^2+a_{3r+3}x^3+\cdots)$, with
$a_{3r+1}, b_{2r+1}, a_{3r+2}, b_{2r+2} \mbox{ and } a_{3r+3}$ not
all zero. The equation of $E_r$ in this chart is $x=0$ and the
numerical data of $E_r$ are $(3r,r+1)$. The zero locus of
$y^3+a_{3r+1}x+b_{2r+1}yx+
a_{3r+2}x^2+b_{2r+2}yx^2+a_{3r+3}x^3+\cdots$ is the strict
transform of $f^{-1}\{0\}$. Remark that it intersects only $E_r$
at this stage.

\vspace{0,2cm} \noindent (b.3.i) If $a_{3r+1} \not= 0$, we obtain
the following after blowing up three more times. \\
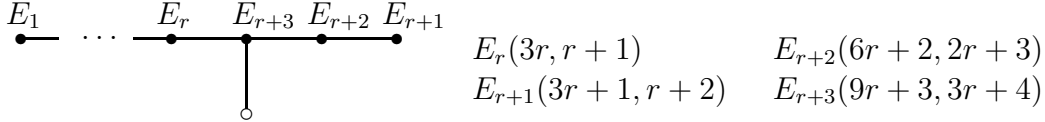
\begin{picture}(150,20)(-5,2)
\put(0,15){\line(1,0){5}} \put(8,15){$\ldots$}
\put(15,15){\line(1,0){35}} \put(0,15){\circle*{1.5}}
\put(20,15){\circle*{1.5}} \put(30,15){\circle*{1.5}}
\put(40,15){\circle*{1.5}} \put(50,15){\circle*{1.5}}
\put(30,15){\line(0,-1){9.25}} \put(30,5){\circle{1.5}}
\put(-2,17){$E_1$} \put(18,17){$E_r$} \put(28,17){$E_{r+3}$}
\put(38,17){$E_{r+2}$} \put(48,17){$E_{r+1}$}
\put(60,7){\shortstack[l]{$E_r(3r,r+1)$ \\ $E_{r+1}(3r+1,r+2)$}}
\put(100,7){\shortstack[l]{$E_{r+2}(6r+2,2r+3)$
\\ $E_{r+3}(9r+3,3r+4)$}}
\end{picture}
The pole $-(3r+4)/(9r+3)$ is in the interval $]-\infty,-1/2]$ if
and only if $r=1$, and in this case the pole is equal to
$-1/2-1/12$.

\vspace{0,2cm} \noindent (b.3.ii) If $a_{3r+1}=0$ and $b_{2r+1}
\not= 0$, calculations give us the following data. \\
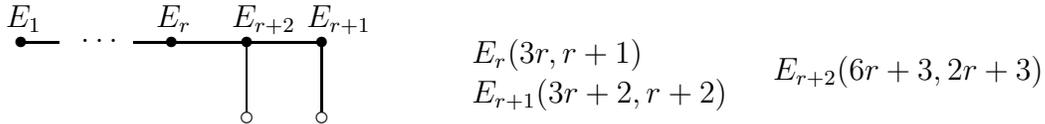
\begin{picture}(150,20)(-5,2)
\put(0,15){\line(1,0){5}} \put(8,15){$\ldots$}
\put(15,15){\line(1,0){25}} \put(0,15){\circle*{1.5}}
\put(20,15){\circle*{1.5}} \put(30,15){\circle*{1.5}}
\put(40,15){\circle*{1.5}} \put(30,15){\line(0,-1){9.25}}
\put(40,15){\line(0,-1){9.25}} \put(30,5){\circle{1.5}}
\put(40,5){\circle{1.5}} \put(-2,17){$E_1$} \put(18,17){$E_r$}
\put(28,17){$E_{r+2}$} \put(38,17){$E_{r+1}$}
\put(60,7){\shortstack[l]{$E_r(3r,r+1)$ \\ $E_{r+1}(3r+2,r+2)$}}
\put(100,10){\shortstack[l]{$E_{r+2}(6r+3,2r+3)$}}
\end{picture}
The pole $-(2r+3)/(6r+3)$ is in the interval $]-\infty,-1/2]$ if
and only if $r=1$, and in this case the pole is equal to
$-1/2-1/18$.

\vspace{0,2cm} \noindent (b.3.iii) If $a_{3r+1}=b_{2r+1}=0$ and
$a_{3r+2} \not= 0$, we get the following. \\
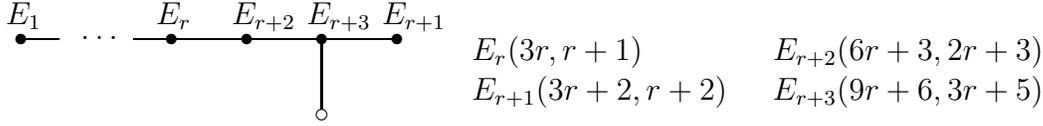
\begin{picture}(150,20)(-5,2)
\put(0,15){\line(1,0){5}} \put(8,15){$\ldots$}
\put(15,15){\line(1,0){35}} \put(0,15){\circle*{1.5}}
\put(20,15){\circle*{1.5}} \put(30,15){\circle*{1.5}}
\put(40,15){\circle*{1.5}} \put(50,15){\circle*{1.5}}
\put(40,15){\line(0,-1){9.25}} \put(40,5){\circle{1.5}}
\put(-2,17){$E_1$} \put(18,17){$E_r$} \put(28,17){$E_{r+2}$}
\put(38,17){$E_{r+3}$} \put(48,17){$E_{r+1}$}
\put(60,7){\shortstack[l]{$E_r(3r,r+1)$ \\ $E_{r+1}(3r+2,r+2)$}}
\put(100,7){\shortstack[l]{$E_{r+2}(6r+3,2r+3)$
\\ $E_{r+3}(9r+6,3r+5)$}}
\end{picture}
The pole $-(3r+5)/(9r+6)$ is in the interval $]-\infty,-1/2]$ if
and only if $r=1$ and in this case the pole is equal to
$-1/2-1/30$.

\vspace{0,2cm} \noindent (b.3.iv) The last case is
$a_{3r+1}=b_{2r+1}=a_{3r+2}=0$ and ($b_{2r+2} \not= 0$ or
$a_{3r+3} \not= 0$). \\ If $y^3+b_{2r+2}yx^2+a_{3r+3}x^3$ is a
product of three distinct linear factors, we get an embedded
resolution after one blowing-up. The numerical data of $E_{r+1}$
are $(3r+3,r+2)$ and $-(r+2)/(3r+3) \notin ]-\infty,-1/2[$. \\ If
$y^3+b_{2r+2}yx^2+a_{3r+3}x^3$ is not a product of three distinct
linear factors, then it is equal to $y^3+xy^2$ after an affine
coordinate transformation that does not change the equation $x=0$
of $E_r$. Let $g_{r+1}$ be the blowing-up at the origin of the
chart we consider. The strict transform of $f^{-1}\{0\}$ only
intersects the exceptional curve $E_{r+1}$, which has numerical
data $(3r+3,r+2)$. Because $-(r+2)/(3r+3) \geq -1/2$ for all $r$,
it follows from (2.4) and (2.6) that $Z_f(s)$ has no pole in
$]-\infty,-1/2[$ different from $-1$ .

\vspace{0,4cm} \noindent \hspace{1cm} (c) Suppose that
$\mbox{mult}(f) \geq 4$. We explained in (2.6) that $Z_f(s)$ has
no pole in $]-\infty,-1/2[$ different from $-1$. $\qquad \Box$

\vspace{0,5cm}

\noindent \textbf{(2.9)} We now present a similar result for the
following generalized zeta functions \cite{DenefLoeser1}. The case
$d=2$ is used in the next section. To $f \in \mathcal{O}_n$ and $d
\in \mathbb{Z}_{>0}$ one associates the local topological zeta
function
\[ Z_f^{(d)}(s) = Z_{\mathrm{top},0,f}^{(d)}(s) := \sum_{I \subset T \atop
\forall i \in I \; : \; d|N_i} \chi(\stackrel{\circ}{E_I} \cap
g^{-1}\{0\}) \prod_{i \in I} \frac{1}{\nu_i+sN_i}.\] For $n,d \in
\mathbb{Z}_{>0}$, we set \[ \mathcal{P}_n^{(d)} := \{ s_0 \mid
\exists f \in \mathcal{O}_n \, : \, Z_{f}^{(d)}(s) \mbox{ has a
pole in } s_0 \}. \] Consequently, $Z_f(s)=Z_f^{(1)}(s)$ and
$\mathcal{P}_n=\mathcal{P}_n^{(1)}$.

\vspace{0,5cm}

\noindent \textbf{(2.10)} Let $E_i$ be an exceptional curve and
let $E_j$, $j \in J$, be the components that intersect $E_i$ in
$V$. Then
\begin{equation}
\label{relatieN} \sum_{j \in J} N_j \equiv 0 \pmod{N_i},
\end{equation}
see e.g. \cite{Loeser} or \cite{Veyscongruences}. Fix $d \in
\mathbb{Z}_{>0}$ and suppose that $d \mid N_i$. Let $J_d \subset
J$ be the subset of indices $j$ satisfying $d \mid N_j$. Suppose
that $\alpha_j:=\nu_j-(\nu_i/N_i)N_j$ is different from $0$ for
all $j \in J_d$. Then the contribution of $E_i$ to the residue
$\mathcal{R}$ of $Z_f^{(d)}(s)$ at the candidate pole $-\nu_i/N_i$
is
\begin{equation} \label{contributiond}
\frac{1}{N_i} \left( \chi(\stackrel{\circ}{E_{\{i\}}}) + \sum_{j
\in J_d} \alpha_j^{-1} \right).
\end{equation}
This contribution is zero if $J$ contains one or two indices.
Indeed, if $J$ contains one element, relation (\ref{relatieN})
implies that $J=J_d$. Therefore, the contribution $\mathcal{R}$ is
the same as in the case $d=1$ and by (2.3) we get $\mathcal{R}=0$.
If $J$ contains two elements, relation (\ref{relatieN}) implies
that $J_d=J$ or $J_d= \emptyset$. If $J_d=J$, we obtain
$\mathcal{R}=0$ analogously as in the previous case. If $J_d=
\emptyset$, we get $\mathcal{R}=0$ because the Euler-Poincar\'{e}
characteristic of a projective line minus two points is zero.

\vspace{0,5cm}

\noindent \textbf{(2.11) Theorem.} \textsl{Let $d \in
\mathbb{Z}_{>1}$. Then}
\[ \mathcal{P}_2^{(d)} \cap \left]-\infty,-\frac{1}{2}\right[ \subset \left\{ \left.
-\frac{1}{2}-\frac{1}{i} \right| i \in \mathbb{Z}_{>1} \right\}.
\]

\vspace{0,2cm} \noindent \emph{Proof.} This follows from the proof
of Theorem 2.8 and from (2.10). $\qquad \Box$

\vspace{0,2cm} \noindent \textsl{Remark.} If one does a lot of
calculations, one can check that
\[ \mathcal{P}_2^{(d)} \cap \left]-\infty,-\frac{1}{2}\right[ = \left\{ \left.
-\frac{1}{2}-\frac{1}{i} \right| i \in \mathbb{Z}_{>2} \mbox{ and
} d|\mbox{lcm}(2,i) \right\} \] if $d \in \mathbb{Z}_{>1}$.
However, we do not need this in the next section.

\section{Surfaces}
In this section, we prove the following theorem.

\vspace{0,2cm}

\noindent \textbf{(3.0) Theorem.} \textsl{We have}
\[ \mathcal{P}_3 \cap ]-\infty,-1[ = \left\{ \left. -1-\frac{1}{i} \right| i \in
\mathbb{Z}_{>1} \right\}. \] \textsl{Moreover, if $f \in
\mathcal{O}_3$ has multiplicity $3$ or more, then $Z_f(s)$ has no
pole less than $-1$.}

\vspace{0,2cm} \noindent \textsl{Remark.} (i) It is a priori not
obvious that the smallest value of $\mathcal{P}_3$ is $-3/2$. This
is in contrast with the fact that it easily follows from lemma 2.5
that the smallest value of $\mathcal{P}_2$ is $-1$.

(ii) In (3.3.9) we give functions $f_k \in \mathcal{O}_3$ of
arbitrary multiplicity such that $Z_{f_k}(s)$ has a pole in $s_k$,
where $(s_k)_k$ is a sequence of real numbers larger than $-1$ and
converging to $-1$. In particular $\mathcal{P}_3 \cap ]-1,-41/42[
\not= \emptyset$, which is in contrast to $\mathcal{T}_3 \cap
]41/42,1[ = \emptyset$.

\subsection{On candidate poles which are not poles}
\noindent \textbf{(3.1.1)} Let $f$ be the germ of a holomorphic
function on a neighbourhood of the origin $0$ in $\mathbb{C}^3$
which satisfies $f(0)=0$ and which is not identically zero. Let
$Y$ be the zero set of $f$. Fix an embedded resolution $g:X_t
\rightarrow X_0 \subset \mathbb{C}^3$ for $Y$ which is an
isomorphism outside the singular locus of $Y$ and which is a
composition $g_1 \circ \cdots \circ g_t$ of blowing-ups $g_i:X_i
\rightarrow X_{i-1}$ with irreducible nonsingular centre $D_{i-1}$
and exceptional variety $E_i^{(0)}$ satisfying for
$i=0,\ldots,t-1$:
\\ \indent (a) the codimension of $D_i$ in $X_i$ is at least 2; \\
\indent (b) $D_i$ is a subset of the strict transform of $Y$ under
$g_1 \circ \cdots \circ g_i$; \\ \indent (c) the union of the
exceptional varieties of $g_1 \circ \cdots \circ g_i$ has only
normal crossings with $D_i$, i.e., for all $P \in D_i$, there are
three surface germs through $P$ which are in normal crossings such
that each exceptional surface germ through $P$ is one of them and
such that the germ of $D_i$ at $P$ is the intersection of some of
them;
\\ \indent (d) the origin $0$ of $\mathbb{C}^3$ is an element of $(g_1 \circ \cdots
\circ g_i)D_i$; and
\\ \indent (e) $D_i$ contains a point in which $(g_1 \circ \cdots
\circ g_i)^{-1}Y$ has not normal crossings. \\ Remark that such a
resolution always exists by Hironaka's theorem \cite{Hironaka}.

\vspace{0,5cm}

\noindent \textbf{(3.1.2)} Fix an exceptional variety $E_i^{(0)}$.
The strict transform $E_i$ of $E_i^{(0)}$ in $X_t$ is obtained by
a finite succession of blowing-ups $h_j$, $j \in
T_e:=\{1,\ldots,m\}$,
\[ E_i^{(0)} \stackrel{h_1}{\longleftarrow} E_i^{(1)}
\stackrel{h_2}{\longleftarrow} \cdots E_i^{(j-1)}
\stackrel{h_j}{\longleftarrow} E_i^{(j)} \cdots
\stackrel{h_{m-1}}{\longleftarrow} E_i^{(m-1)}
\stackrel{h_m}{\longleftarrow} E_i^{(m)}=E_i \] with centre
$P_{j-1} \in E_i^{(j-1)}$ and exceptional curve $C_j^{(j)}$. The
irreducible components of the intersection of $E_i^{(0)}$ with
irreducible components of $(g_1 \circ \cdots \circ g_i)^{-1}Y$
different from $E_i^{(0)}$ are denoted by $C_j^{(0)}$, $j \in
T_s$. The strict transform of $C_j^{(k)}$ in $E_i^{(l)}$ is
denoted (whenever this makes sense) by $C_j^{(l)}$ and we set
$C_j=C_j^{(m)}$. Remark that $h:=h_1 \circ \cdots \circ h_m$ is an
embedded resolution of $\cup_{j \in T_s}C_j^{(0)}$. For each $j
\in T:=T_s \cup T_e$ the curve $C_j$ is an irreducible component
of the intersection of $E_i$ with exactly one other component of
$g^{-1}Y$. Let this component have numerical data $(N_k,\nu_k)$
and set $\alpha_j=\nu_k-(\nu_i/N_i)N_k$.

\vspace{0,5cm}

\noindent \textbf{(3.1.3)} Suppose that $E_i^{(0)} \subset (g_1
\circ \cdots \circ g_i) ^{-1}\{0\}$ and that $\alpha_j \not= 0$
for every $j \in T$. The contribution $\mathcal{R}$ of $E_i$ to
the residue of $Z_f(s)$ at the candidate pole $-\nu_i/N_i$ is
\begin{equation}\label{contributionsurface}
\frac{1}{N_i}\left(\sum_{I \subset T} \chi(\stackrel{\circ}{C_I})
\prod_{j \in I} \alpha_j^{-1}\right),
\end{equation}
where $\stackrel{\circ}{C_I}$ denotes the subset $(\cap_{j \in I}
C_j) \setminus (\cup_{j \not\in I} C_j)$ of $E_i$. Remark that
$\stackrel{\circ}{C_{\emptyset}}=E_i \setminus (\cup_{j \in T}
C_i)$. We now state some relations between the $\alpha_i$, which
will allow us to prove that this contribution is identically zero
(i.e., zero for any value of the alphas) for a lot of intersection
configurations on $E_i^{(0)}$.

\vspace{0,5cm}

\noindent \textbf{(3.1.4)} To the creation of $E_i^{(0)} \subset
(g_1 \circ \cdots \circ g_i)^{-1}\{0\}$ in the resolution process,
we associate the relation
\begin{equation}\label{relcreation}
\sum_{j \in T_s} d_j(\alpha_j-1)+3-\dim D_{i-1}=0,
\end{equation}
where $d_i$, $i \in T_s$, is the degree of the intersection cycle
$C_i^{(0)} \cdot F$ on $F$ for a general fibre $F$ of
$g_i|_{E_i^{(0)}}:E_i^{(0)} \rightarrow D_{i-1}$ over a point of
$D_{i-1}$. In particular, when $D_{i-1}$ is a point, we have that
$E_i^{(0)} \cong \mathbb{P}^2$ and that $d_i$ is just the degree
of the curve $C_i^{(0)}$. To the blowing-up $h_j$ we associate the
relation
\begin{equation}\label{relblowingup}
\alpha_j= \sum_{ k \in T_s \cup \{1,\ldots,j-1\} }
\mu_k(\alpha_k-1)+2,
\end{equation}
where $\mu_k$, $k \in T_s \cup \{1,\ldots,j-1\}$, is the
multiplicity of $P_{j-1}$ on $C_k^{(j-1)}$. See
\cite{Veysrelations} for more general statements in arbitrary
dimension and proofs.

\vspace{0,5cm}

\noindent \textbf{(3.1.5)} Now we proceed in the same way as in
\cite{Veysconfigurations} for Igusa's $p$-adic zeta function. One
easily verifies that the number (\ref{contributionsurface}) does
not change when we do an extra blowing-up $h_{m+1}$ at a point
$P_m \in E_i^{(m)}$ and associate to the new exceptional curve a
number $\alpha$ using (\ref{relblowingup}). Because of this
observation, one can compute $\mathcal{R}$ if one has the curves
$C_j^{(0)}$, $j \in T_s$, on $E_i^{(0)}$ together with the
associated values $\alpha_j$ as follows. Compute the
\textsl{minimal} embedded resolution of $\cup_{j \in T_s}
C_j^{(0)}$ and compute the alpha associated to an exceptional
curve using (\ref{relblowingup}). By putting these data in
(\ref{contributionsurface}), we get $\mathcal{R}$.

\vspace{0,5cm}

\noindent \textbf{(3.1.6)} \textsl{Example.} Suppose that
$E_i^{(0)}$ is the exceptional variety of a blowing-up at a point
and suppose that the intersection configuration on $E_i^{(0)}$
consists of three projective lines $C_j^{(0)}$, $j \in
T_s:=\{a,b,c\}$, all passing through the same point $P$. Suppose
that $\alpha_j \not= 0$ for all $j \in T$. The minimal embedded
resolution $l:W \rightarrow E_i^{(0)}$ is the blowing-up at $P$.
By abuse of notation, we denote the exceptional curve by $C_1$ and
the strict transform of $C_j^{(0)}$, $j \in T_s$, by $C_j$.

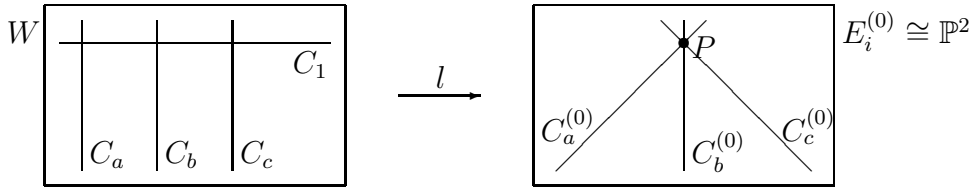
\begin{picture}(145,34)(15,0)
\put(20,5){\line(1,0){40}} \put(20,29){\line(1,0){40}}
\put(20,5){\line(0,1){24}} \put(60,5){\line(0,1){24}}
\put(85,5){\line(1,0){40}} \put(85,29){\line(1,0){40}}
\put(85,5){\line(0,1){24}} \put(125,5){\line(0,1){24}}
\put(25,7){\line(0,1){20}} \put(35,7){\line(0,1){20}}
\put(45,7){\line(0,1){20}} \put(22,24){\line(1,0){36}}
\put(88,7){\line(1,1){20}} \put(122,7){\line(-1,1){20}}
\put(105,7){\line(0,1){20}} \put(67,17){\vector(1,0){11}}
\put(15,24){\shortstack[l]{$W$}}
\put(126,24){\shortstack[l]{$E_i^{(0)} \cong \mathbb{P}^2$}}
\put(72,18){\shortstack[l]{$l$}} \put(26,8){\shortstack[l]{$C_a$}}
\put(36,8){\shortstack[l]{$C_b$}}
\put(46,8){\shortstack[l]{$C_c$}}
\put(53,20){\shortstack[l]{$C_1$}}
\put(106,22){\shortstack[l]{$P$}} \put(105,24){\circle*{1.5}}
\put(86,11){\shortstack[l]{$C_a^{(0)}$}}
\put(106,8){\shortstack[l]{$C_b^{(0)}$}}
\put(118,11){\shortstack[l]{$C_c^{(0)}$}}
\end{picture}
By relations (\ref{relcreation}) and (\ref{relblowingup}) we have
$\alpha_a+\alpha_b+\alpha_c=0$ and
$\alpha_1=\alpha_a+\alpha_b+\alpha_c-1=-1$ respectively. Now we
can calculate the contribution $\mathcal{R}$ of the strict
transform of $E_i^{(0)}$ in $X_t$ to the residue of $Z_f(s)$ at
the candidate pole $-\nu_i/N_i$:
\begin{eqnarray*}
\mathcal{R} & = & \frac{1}{N_i} \left( \sum_{I \subset T}
\chi(\stackrel{\circ}{C_I})\prod_{j \in I} \alpha_j^{-1} \right)
\\ & = & \frac{1}{N_i}
\left(-1-\frac{1}{\alpha_1} +\frac{1}{\alpha_a}
+\frac{1}{\alpha_b} +\frac{1}{\alpha_c}
+\frac{1}{\alpha_1\alpha_a} +\frac{1}{\alpha_1\alpha_b}
+\frac{1}{\alpha_1\alpha_c} \right) \\  & = & 0.
\end{eqnarray*}
We stress that $\mathcal{R}$ is zero for any possible value of
$\alpha_a$, $\alpha_b$ and $\alpha_c$.

\subsection{Multiplicity 2}
\noindent \textbf{(3.2.1)} Let $f$ be the germ of a holomorphic
function on a neighbourhood of the origin $0$ in $\mathbb{C}^n$
which satisfies $f(0)=0$, and let $F$ be the germ of the
holomorphic function $f+x^2_{n+1}$ on a neighbourhood of the
origin $0$ in $\mathbb{C}^{n+1}$. Then the following equality is
obtained in \cite{Artal}, see also the Thom-Sebastiani principle
in \cite{DenefLoeser3}:
\begin{eqnarray*}
Z_F(s) & = & \frac{1}{2s+1} + \frac{s(2s+3)}{2(s+1)(2s+1)}Z_f
\left( s+\frac{1}{2} \right) - \frac{3s}{2(s+1)}Z_f^{(2)} \left(
s+\frac{1}{2} \right).
\end{eqnarray*}

\vspace{0,5cm}

\noindent \textbf{(3.2.2) Proposition.} \textsl{The set}
\[ \left\{ s_0 \mid \exists f \in \mathcal{O}_3 \mbox{ : mult}(f)=2 \textsl{\mbox{
and }} Z_f(s) \textsl{\mbox{ has a pole in }} s_0 \right\} \cap
\left] -\infty,-1 \right[ \] \textsl{is equal to}
\[ \left\{ \left. -1-\frac{1}{i} \right| i \in \mathbb{Z}_{>1} \right\}. \]

\vspace{0,2cm} \noindent \emph{Proof.} Let $f$ be an element of
$\mathcal{O}_3$ with  multiplicity $2$. Up to an affine coordinate
transformation, the part of degree two in the Taylor series of $f$
is equal to $x^2$, $x^2+y^2$ or $x^2+y^2+z^2$. Using (2.7), we may
suppose that $f$ is of the form $x^2+g(y,z)$ with $g(y,z) \in
\mathcal{O}_2$. The formula in (3.2.1) and the result for curves
imply that every pole of $Z_f(s)$ less than $-1$ is of the form
$-1-1/i$, $i \in \mathbb{Z}_{>1}$. For the other inclusion, we
remark that the poles of the local topological zeta function
associated to $x^2+y^2+z^i$, $i \geq 2$, are $-1-1/i$ and $-1$.
$\qquad \Box$

\vspace{0,5cm}

\noindent \textbf{(3.2.3)} Our next goal is to give a sequence of
poles larger than $-1$ and converging to $-1$. Keeping in mind the
formula in (3.2.1), we try to find functions $f_k \in
\mathcal{O}_2$ such that $Z_{f_k}(s)$ has a pole in $s_k$, where
$(s_k)_k$ is a sequence of real numbers larger than $-1/2$ and
converging to $-1/2$. Set $f_k=x^3y^2+x^k$ for $k \geq 5$.
\\ We obtain the following equalities after some calculations:
\begin{eqnarray*}
 Z_{f_{2r+4}}(s)=\frac{3s^2+2rs+8s+2r+3}{(4rs+8s+2r+3)(3s+1)(s+1)} \; , &
  \! \! \! \! \! \! & Z_{f_{2r+4}}^{(2)}(s)=\frac{1}{4rs+8s+2r+3} \; , \\
Z_{f_{2r+3}}(s)=\frac{3s^2-rs-2s-r-1}{(2rs+3s+r+1)(3s+1)(s+1)} \;
, & \! \! \! \! \! \! & Z_{f_{2r+3}}^{(2)}(s)=0.
\end{eqnarray*}
Now we use the formula in (3.2.1) to calculate the local
topological zeta function of $F_k:=f_k+z^2$.  We obtain for even
and odd $k$ that
\[Z_{F_k}(s)=\frac{(6k-6)s^2+(15k-5)s+10k-5}{(6s+5)(s+1)(2ks+2k-1)}.\]
Finally, we make the substitution $s=-(2k-1)/(2k)$ in the
numerator in order to check that this value, which converges to
$-1$ if $k$ goes to infinity, is a pole. We obtain
\[ \frac{(k-1)(k-3)(2k-1)}{2k^2}. \]
This value never becomes zero because $k \geq 5$. Consequently,
$-(2k-1)/(2k)$ is always a pole of $Z_{F_k}(s)$.

\vspace{0,2cm} \noindent \textsl{Remark.} In particular we obtain
that $\mathcal{P}_3 \cap ]-1,-41/42[ \not= \emptyset$, which is in
contrast to $\mathcal{T}_3 \cap ]41/42,1[ = \emptyset$.

\subsection{Multiplicity larger than 2}
\noindent \textbf{(3.3.1)} Let $f$ be the germ of a holomorphic
function on a neighbourhood of the origin $0$ in $\mathbb{C}^3$
which satisfies $f(0)=0$ and which is not identically zero. Let
$Y$ be the zero set of $f$. Fix an embedded resolution $g$ for $Y$
which is a composition of blowing-ups $g_{ij}:X_i \rightarrow X_j$
with irreducible nonsingular centre $D_j$ and exceptional surface
$E_i$ as in (3.1.1). Denote the irreducible components of $Y$ by
$E_i$, $i \in T_s$. The strict transform of a variety $E_i$ by a
succession of blowing-ups will be denoted in the same way. The
numerical data of $E_i$ are $(N_i,\nu_i)$.

\vspace{0,5cm}

\noindent \textbf{(3.3.2)} The following table gives the numerical
data of $E_i$. In the columns, the dimension of $D_j$ is kept
fixed. In the rows, the number of exceptional surfaces through
$D_j$ is kept fixed. So $E_k$, $E_l$ and $E_m$ represent
exceptional surfaces that contain $D_j$. The multiplicity of $D_j$
on the strict transform of $Y$ is denoted by $\mu_{D_j}$.
\begin{picture}(145,20)(0,0)
\put(0,-12){\line(0,1){21}} \put(31,14){\line(1,0){113}}
\begin{tabular}{c|c|c|}
  & $D_j$ is a point $P$ & $D_j$ is a curve $L$ \\ \hline
 / & $(\mu_P,3)$ & $(\mu_L,2)$ \\ \hline
 $E_k$ & $(N_k+\mu_P,\nu_k+2)$ & $(N_k+\mu_L,\nu_k+1)$ \\ \hline
 $E_k$ and $E_l$ & $(N_k+N_l+\mu_P,\nu_k+\nu_l+1)$ &
 $(N_k+N_l+\mu_L,\nu_k+\nu_l)$ \\ \hline
 $E_k$, $E_l$ and $E_m$ & $(N_k+N_l+N_m+\mu_P,\nu_k+\nu_l+\nu_m)$
 & / \\ \hline
\end{tabular}
\end{picture}

\vspace{1,5cm}

\noindent \textbf{(3.3.3) Lemma.} \textsl{Suppose that} mult$(f)
\geq 3$. \textsl{If there is no exceptional surface through $D_j$,
then $-\nu_i/N_i \geq -1$.}

\vspace{0,2cm} \noindent \emph{Proof.} The case that the centre
$D_j$ is a point $P$ through which no exceptional surface passes
can only occur in the first blowing-up because of condition (d) in
(3.1.1) and because the inverse image of $0$ in $X_j$ is contained
in the union of the exceptional surfaces in $X_j$. Since
$\mbox{mult}(f) \geq 3$, we have in this case
$-\nu_i/N_i=-3/\mu_P=-3/\mbox{mult}(f) \geq -1$. \\ If the centre
$D_j$ is a curve $L$ contained in no exceptional surface, then
$\mu_L \geq 2$ because our embedded resolution is an isomorphism
outside the singular locus of $Y$. Consequently, we get in this
case $-\nu_i/N_i=-2/\mu_L \geq -1$. $\qquad \Box$

\vspace{0,5cm}

\noindent \textbf{(3.3.4)} Suppose that $D_j$ is contained in at
least one exceptional surface and that the candidate poles
associated to the exceptional surfaces that pass through $D_j$ are
larger than or equal to $-1$. Then the table in (3.3.2) implies
that also $-\nu_i/N_i \geq -1$, unless $D_j$ is a nonsingular
point $P$ of the strict transform of $Y$ through which only one
exceptional surface $E_0$ passes and $-\nu_0/N_0=-1$. Suppose that
we are in this situation. Denote the unique irreducible component
of the strict transform of $Y$  which passes through $P$ by $E_a$.
Consider now a small enough neighbourhood $Z_0$ of $P$ on which
$E_a$ is nonsingular such that, if we restrict the blowing-ups
$g_{ij}$ to the inverse image of $Z_0$, we get an embedded
resolution $h=h_1 \circ \cdots \circ h_s$ for the germ of $E_a
\cup E_0$ at $P$ which is a composition of blowing-ups $h_i:Z_i
\rightarrow Z_{i-1}$, $i \in \{1,\ldots,s\}$, with irreducible
nonsingular centre $D_{i-1}':=D_{i-1} \cap Z_{i-1}$ and
exceptional surface $E_i':=E_i \cap Z_i$ satisfying for
$i=0,\ldots,s-1$:

\vspace{0,1cm} \indent (a) the codimension of $D_i'$ in $Z_i$ is
at least $2$;
\\ \indent (b) $D_i'$ is a subset of $E_a':=E_a \cap Z_i$;
\\ \indent (c) $\cup_{l \in \{0,1,\ldots,i\}} E_l'$ has only normal
crossings with $D_i'$, where $E_0':=E_0 \cap Z_0$; \\ \indent (d)
the image of $D_i'$ under $h_1 \circ \cdots \circ h_i$ contains
$P$; and \\ \indent (e) if $D_i=D_i'$, then $D_i$ contains a point
where there are no normal crossings.

\vspace{0,2cm} \noindent Remark that it can happen that $g_{ij}$
is an isomorphism on the inverse image of $Z_0$. Because we did
not specify the indices in (3.3.1), we were able to get a nice
notation here. Remark also that $D_i=D_i'$ if $D_i$ is a point.
From now on, we study the resolution $h:Z_s \rightarrow Z_0$ for
the germ of $E_a \cup E_0$ at $P$.

\vspace{0,5cm}

\noindent \textbf{(3.3.5) Lemma.} \textsl{If $D_i=D_i'$, then
$D_i$ is a subset of $E_0'$.}

\vspace{0,2cm} \noindent \emph{Proof.} Remark that $D_i$ has to
lie in an exceptional surface because $E_a'$ is nonsingular and
because an embedded resolution is an isomorphism outside the
singular locus of $Y$.

\noindent \hspace{0,5cm} First we consider the case that
$D_i=D_i'$ is a point contained in exceptional surfaces different
from $E_0'$ and in the surface $E_a'$. The union of these surfaces
has normal crossings at $D_i$ because $E_a'$, considered as a
subset of $Z_0$, is nonsingular. This is in contradiction with
(e). Remark that it can thus not happen that $E_a'$ and three
exceptional surfaces different from $E_0'$ have a point in common.

\noindent \hspace{0,5cm} The case that $D_i=D_i'$ is a curve
contained in exactly two exceptional surfaces different from
$E_0'$ and in the surface $E_a'$ cannot occur because $E_a'$ is a
nonsingular subset of $Z_0$ and therefore these three surfaces
should have normal crossings.

\noindent \hspace{0,5cm} Finally we study the case that $D_i=D_i'$
is a curve contained in one exceptional surface $E_j'$ different
from $E_0'$ and in $E_a'$. Condition (c) implies that every point
of $D_i$ is contained in at most one exceptional surface different
from $E_j'$. Moreover, such an exceptional surface has to be
transversal to $D_i$. This implies that there are normal crossings
at every point of $D_i$, which is in contradiction with (e).
Therefore, this case cannot occur. $\qquad \Box$

\vspace{0,5cm}

\noindent \textbf{(3.3.6) Lemma.} \textsl{Suppose that} mult$(f)
\geq 3$. \textsl{Then we have} $\nu_i \leq N_i+1$ \textsl{for
every exceptional surface} $E_i$, $i \in \{1,\ldots,s\}$.
\textsl{Moreover,} $\nu_i=N_i+1$ \textsl{if and only if} $D_{i-1}$
\textsl{is a point and the numerical data of every exceptional
surface} $E_j$ \textsl{different from} $E_0$ \textsl{and through}
$D_{i-1}$ \textsl{satisfy} $\nu_j=N_j+1$.

\vspace{0,2cm} \noindent \emph{Proof.} The proof is by induction
on $i$. Since $\nu_0=N_0$, we have that $\nu_1=N_1+1$. Suppose now
that $\nu_j \leq N_j+1$ for every exceptional surface $E_j$
through $D_{i-1}$.

\noindent \hspace{0,5cm} \textsl{Case 1: $D_{i-1}$ is a point.} We
obtain from (3.3.5) that $D_{i-1}$ is a subset of $E_0'$. Because
$\nu_0=N_0$ and because every other exceptional surface $E_j$
through $D_{i-1}$ satisfies $\nu_j \leq N_j+1$, the table of
(3.3.2) gives us that $\nu_{i} \leq N_{i}+1$.

\noindent \hspace{0,5cm} \textsl{Case 2: $D_{i-1}$ is a curve.} If
$D_{i-1} \not= D_{i-1}'$, then $D_{i-1}' \not\subset (h_1 \circ
\cdots \circ h_{i-1})^{-1}P$ and therefore we get as in the
beginning of (3.3.4) that $-\nu_i/N_i \geq -1$. If
$D_{i-1}=D_{i-1}'$, one computes from (3.3.2) and the previous
lemma that $-\nu_i/N_i \geq -1$.

\noindent \hspace{0,5cm} We have now proved the first part of the
lemma. Using this first part and the table of (3.3.2), we get the
second part. $\qquad \Box$

\vspace{0,5cm}

\noindent \textbf{(3.3.7) Lemma.} \textsl{If} mult$(f) \geq 3$
\textsl{and if the numerical data of} $E_i$ \textsl{satisfy}
$\nu_i=N_i+1$, \textsl{then} $-\nu_i/N_i \not= -\nu_j/N_j$
\textsl{for every exceptional surface} $E_j$ \textsl{that
intersects} $E_i$ \textsl{at some stage of the resolution
process.}

\vspace{0,2cm} \noindent \emph{Proof.} Let $E_j$ be an exceptional
surface that intersects $E_i$ at some stage of the resolution
process. If $E_j$ is created before $E_i$, then $E_j$ contains the
point $D_{i-1}$. Otherwise, $E_j$ is created by a blowing-up at a
point of $E_i$ or by a blowing-up along a curve.

\noindent \hspace{0,5cm} If $E_j$ is created by a blowing-up along
a curve, then $-\nu_j/N_j \geq -1$, and consequently $-\nu_i/N_i
\not= -\nu_j/N_j$. Now we consider the case that $E_j$ contains
the point $D_{i-1}$. There is no problem if $\nu_j \leq N_j$.
Consequently, suppose that $\nu_j=N_j+1$. From the table in
(3.3.2), we get $N_j < N_i$. Therefore, $-\nu_i/N_i=-(N_i+1)/N_i >
-(N_j+1)/N_j=-\nu_j/N_j$. The case that $E_j$ is created by a
blowing-up at a point of $E_i$ is treated analogously. $\qquad
\Box$

\vspace{0,5cm}

\noindent \textbf{(3.3.8) Proposition.} \textsl{If} mult$(f) \geq
3$, \textsl{then no pole of $Z_f(s)$ is less than $-1$.}

\vspace{0,2cm} \noindent \emph{Proof.} Suppose that mult$(f) \geq
3$.

\noindent \hspace{0,5cm} We have only to consider exceptional
surfaces with a candidate pole less than $-1$. Recall from (3.3.6)
that $-\nu_i/N_i <-1$ if and only if $D_{i-1}$ is a point and all
exceptional surfaces through the point $D_{i-1}$ different from
$E_0$ have a candidate pole less than $-1$. We will determine all
possible intersection configurations on such surfaces just after
their creation.

\noindent \hspace{0,5cm} If $-\nu_i/N_i \geq -1$ and
$-\nu_{i+1}/N_{i+1} < -1$, then the blowing-ups along $D_{i-1}$
and $D_i$ commute with each other. Therefore, we may assume that
there is a $k$ (larger than zero because $-\nu_1/N_1<-1$) such
that $-\nu_i/N_i < -1$ for $1 \leq i \leq k$ and $-\nu_i/N_i \geq
-1$ for $k<i \leq s$.

\noindent \hspace{0,5cm} The intersection configuration on $E_1$
consists of one projective line, which is the intersection with
$E_0$ and $E_a$. The points of $Z_1$ in which we do not have
normal crossings and which lie above $P$ are those on this
projective line. This implies the following statement for $i=2$.
\begin{quote}
If $Q$ is a point of $Z_{i-1}$, $i \in \{2,\ldots,k\}$, in which
we do not have normal crossings and which lies above $P$ (so
consequently $Q$ is a point of $E_0$, of one or two other
exceptional surfaces and of $E_a$), then there exists an
exceptional surface $E_l$ through $Q$ with the property $E_0 \cap
E_l = E_a \cap E_l$.
\end{quote}
\vspace{-0,5cm}
\begin{picture}(150,0)(0,0)
\put(138,12){($*$)}
\end{picture}
We prove this statement by induction on $i$. Suppose that it is
true for $i=j \in \{2, \ldots ,k-1\}$. We give the proof for
$i=j+1$. The statement follows from the induction hypothesis for
points not on $E_j$, because a blowing-up is an isomorphism
outside the exceptional surface. So we prove it for points on
$E_j$. By the induction hypothesis applied to the point $D_{j-1}$,
we obtain that there exists an exceptional surface $E_l$ through
$D_{j-1}$ such that $E_0 \cap E_l= E_a \cap E_l$ in $Z_{j-1}$. But
then $E_a \cap E_l=E_0 \cap E_l$ in $Z_j$, which solves the
problem for the point $E_0 \cap E_l \cap E_j$. There are other
points on $E_j$ in which we do not have normal crossings if and
only if $E_a$ is tangent to $E_0$ in $D_{j-1}$. In this case, the
points in which we do not have normal crossings are the points of
$E_0 \cap E_j$. Because $E_0 \cap E_j = E_a \cap E_j$, we are
done.

\noindent \hspace{0,5cm} Because the centre of a blowing-up
satisfies the conditions of the statement, we obtain that the
possible intersection configurations are the following
configurations of lines in $\mathbb{P}^2$: (i) one line, (ii) two
lines, (iii) three lines through one point, (iv) three lines in
general position and (v) three lines through one point and a
fourth line not through that point.

\noindent \hspace{0,5cm} For all these configurations, we can
calculate as in (3.1.6) that the contribution to the residue is
$0$. The second author did this already in
\cite{Veysconfigurations} for Igusa's $p$-adic zeta function. The
point is that ($*$) excludes the configuration consisting of four
lines in general position, for which this contribution is not
zero. Remark also that we need here that the alphas are not zero,
a fact we proved in (3.3.7). $\qquad \Box$

\vspace{0,5cm}

\noindent \textbf{(3.3.9)} In (3.2.3), we found functions $f_k \in
\mathcal{O}_3$ of multiplicity $2$ such that $Z_{f_k}(s)$ has a
pole in $s_k$, where $(s_k)_k$ is a sequence of real numbers
larger than $-1$ and converging to $-1$. Here we construct for
\textsl{every} $n \geq 0$ functions $f_k \in \mathcal{O}_3$ of
multiplicity $n+2$ with this property. We use the formula obtained
by Denef and Loeser in \cite[Th\'eor\`eme~5.3]{DenefLoeser1},
which expresses the local topological zeta function of a
non-degenerated polynomial in terms of its Newton polyhedron. Fix
$n \geq 0$ and set $f_k=x^nz^2+x^{3+n}y^2+x^k$ for $k \geq n+4$.
Then
\[ Z_{f_k}(s)=\frac{\begin{array}{c}(-2n^2-6n)s^3+(n^2+3kn-4n+6k-6)s^2 \\ +(-4n^2
+4kn-7n+15k-5)s-10n+10k-5\end{array}}{(6s+2ns+5)(s+1)(2ks+2k-2n-1)(ns+1)}.
\] Consequently, $-(2k-2n-1)/(2k)$ is a pole if and only if it is
not a zero of the numerator. So we make the substitution
$s=-(2k-2n-1)/(2k)$ in the numerator and obtain
\[ \frac{(k-1-2n)(k-n-3)(2k-2n-1)(2n^2-2kn+n+2k)}{4k^3}. \]
Because $k \geq n+4$, this is zero if and only if $k=1+2n$. Thus
we have found for any multiplicity larger than one a sequence with
the desired property.

\section{Other zeta functions}
\noindent \textbf{(4.1)} Denef and Loeser associate in
\cite{DenefLoeser2} to a polynomial its motivic zeta function,
which is a much finer invariant than its topological zeta
function. Instead of the usual topological Euler-Poincar\'{e}
characteristic, it involves the so-called universal Euler
characteristic of an algebraic variety, i.e., its class in the
Grothendieck ring.

We recall this notion. The Grothendieck ring
$K_0(\mathrm{Var}_{\mathbb{C}})$ of complex algebraic varieties is
the free abelian group generated by the symbols $[V]$, where $V$
is a variety, subject to the relations $[V]=[V']$, if $V$ is
isomorphic to $V'$, and $[V]=[V \setminus W]+[W]$, if $W$ is
closed in $V$. Its ring structure is given by $[V] \cdot [W] := [V
\times W]$. We set $\mathbb{L}:=[\mathbb{A}^1_{\mathbb{C}}]$ and
denote by $\mathcal{M}$ the localization of
$K_0(\mathrm{Var}_{\mathbb{C}})$ with respect to $\mathbb{L}$.

\vspace{0,5cm}

\noindent \textbf{(4.2)} In \cite{DenefLoeser2} the motivic zeta
function is more generally defined for a regular function $f$ on a
smooth algebraic variety $X$, with respect to a subvariety $W$ of
$X$; we refer to \cite[section 2]{DenefLoeser2} for this
definition. One easily verifies that the construction is still
valid for a germ $f$ of a holomorphic function at $0 \in
\mathbb{C}^n$ when $W=\{0\}$; we denote this (local) motivic zeta
function by $Z_{\mathrm{mot},0,f}(s)$. Then, with the notation of
(1.1), the formula of \cite[Theorem 2.2.1]{DenefLoeser2} yields
that \[ Z_{\mathrm{mot},0,f}(s)=\mathbb{L}^{-n} \sum_{I \subset T}
[\stackrel{\circ}{E_I} \cap g^{-1}\{0\}] \prod_{i \in I}
\frac{\mathbb{L}-1}{\mathbb{L}^{\nu_i+sN_i}-1}. \] Here
$\mathbb{L}^{-s}$ should be considered as a variable, and this
expression lives in a localization of the polynomial ring
$\mathcal{M}[\mathbb{L}^{-s}]$.

\vspace{0,5cm}

\noindent \textbf{(4.3)} The motivic zeta function
$Z_{\mathrm{mot},0,f}(s)$ specializes to $Z_{\mathrm{top},0,f}(s)$
\cite[subsection 2.3]{DenefLoeser2}, but also to various
`intermediate level' zeta functions. An important one uses Hodge
polynomials. Recall that the Hodge polynomial of a complex
algebraic variety $V$ is \[ H(V)=H(V,u,v):= \sum_{p,q} \left(
\sum_{i \geq 0} (-1)^i h^{p,q}\left(H_c^i(V,\mathbb{C})\right)
\right)u^pv^q \in \mathbb{Z}[u,v],\] where
$h^{p,q}\left(H^i_c(V,\mathbb{C})\right)$ is the rank of the
$(p,q)$-Hodge component of the $i$-th cohomology group with
compact support of $V$. The zeta function of $f$ on this level is
\[ Z_{\mathrm{Hod},0,f}(s)=(uv)^{-n} \sum_{I \subset T}
H\left(\stackrel{\circ}{E_I} \cap g^{-1}\{0\}\right) \prod_{i \in
I} \frac{uv-1}{(uv)^{\nu_i+sN_i}-1}; \] here $(uv)^{-s}$ is a
variable, and this zeta function lives e.g. in the field of
rational functions in $(uv)^{-s}$ over $\mathbb{Q}(u,v)$.

\vspace{0,5cm}

\noindent \textbf{(4.4)} As in \cite{RodriguesVeys} we define the
poles of $Z_{\mathrm{Hod},0,f}(s)$ to be the real numbers $s_0$
such that $(uv)^{-s_0}$ is a pole of $Z_{\mathrm{Hod},0,f}(s)$,
considered as rational function in $(uv)^{-s}$. Then we have the
following.

Theorems 2.8 and 3.0 are still valid with
$Z_f(s)=Z_{\mathrm{top},0,f}(s)$ replaced by
$Z_{\mathrm{Hod},0,f}(s)$ and $\mathcal{P}_n=\{s_0 \mid \exists f
\in \mathcal{O}_n \; : \; Z_{\mathrm{Hod},0,f}(s)$ has a pole in
$s_0\}$. The proofs are the same as before; they essentially just
use the `geometry' of a resolution.

A good definition of poles of $Z_{\mathrm{mot},0,f}(s)$ is not
immediately clear, due to the fact that $\mathcal{M}$ could have
zero divisors (at present this is an open question). Using the
definition of \cite{RodriguesVeys} for real poles, Theorems 2.8
and 3.0 are also valid for $Z_{\mathrm{mot},0,f}(s)$.

\vspace{0,5cm}

\noindent \textbf{(4.5)} One could and should also wonder whether
there are analogous results for Igusa's $p$-adic zeta function.
This problem is studied in a next paper \cite{Segers}.

\footnotesize{

\noindent \textsc{K.U.Leuven, Departement Wiskunde,
Celestijnenlaan 200B, B-3001 Leuven, Belgium,} \\ \textsl{E-mail:}
dirk.segers@wis.kuleuven.ac.be, wim.veys@wis.kuleuven.ac.be,\\
http://www.wis.kuleuven.ac.be/algebra/veys.htm}

\end{document}